\documentclass[10pt, titlepage]{amsart}

\usepackage{ae} 
\usepackage[T1]{fontenc}
\usepackage[cp1250]{inputenc}
\usepackage{amsmath}
\usepackage{amssymb, amsfonts,amscd,verbatim}
\usepackage[dvips]{graphicx}
\usepackage{latexsym}
\usepackage{indentfirst}
\usepackage{latexsym}

\usepackage{graphicx}
\usepackage{verbatim}   
\usepackage{color}      
\usepackage{subfigure}  

\usepackage{amsmath}    

\theoremstyle{plain}
\newtheorem{Prop}{Proposition}[section]
\newtheorem{Thm}[Prop]{Theorem}
\newtheorem{Cor}[Prop]{Corollary}

\newtheorem{Lem}[Prop]{Lemma}

\theoremstyle{definition}
\newtheorem{Def}[Prop]{Definition}

\theoremstyle{remark}

\def\NN{{\mathbf N}}

\def\UU{{\mathcal U}}
\def\VV{{\mathcal V}}
\def\WW{{\mathcal W}}

\def\asdim{\mathrm{asdim}}
\def\dim{\mathrm{dim}}

\def\dokaz{{\bf Proof. }}
\def\edokaz{\hfill $\blacksquare$}

\numberwithin{equation}{section}

\input pstricks.tex
\input xy
\xyoption{all}

\begin{document}
\title[
Asymptotic dimension of coarse spaces via maps to simplicial complexes
]%
   {Asymptotic dimension of coarse spaces via maps to simplicial complexes}

\author{M.~Cencelj}
\address{IMFM, Pedago\v ska fakulteta, Jadranska ulica 19,
SI-1111 Ljubljana,
Slovenija }
\email{matija.cencelj@guest.arnes.si}

\author{J.~Dydak}
\address{University of Tennessee, Knoxville, TN 37996, USA}
\email{jdydak@utk.edu}

\author{A.~Vavpeti\v c}
\address{Fakulteta za Matematiko in Fiziko,
Univerza v Ljubljani,
Jadranska ulica 19,
SI-1111 Ljubljana,
Slovenija }
\email{ales.vavpetic@fmf.uni-lj.si}

\date{ \today
}
\keywords{asymptotic dimension, coarse geometry, Lipschitz maps, Property A}

\subjclass[2000]{Primary 54F45; Secondary 55M10}

\thanks{This research was supported by the Slovenian Research
Agency grants P1-0292-0101, J1-6721-0101, J1-5435-0101}

\begin{abstract}

It is well-known that a paracompact space $X$ is of covering dimension at most $n$
if and only if any map $f\colon X\to K$ from $X$ to a simplicial complex $K$ can be pushed into its $n$-skeleton
$K^{(n)}$. We use the same idea to characterize asymptotic dimension in the coarse category
of arbitrary coarse spaces. Continuity of the map $f$ is replaced by variation of $f$ on elements
of a uniformly bounded cover. The same way one can generalize Property A of G.Yu to arbitrary coarse spaces.
\end{abstract}

\maketitle

\medskip
\medskip
\tableofcontents
\section{Introduction}

It is well-known (see \cite{Dyd}) that the covering dimension $\dim(X)$
of a paracompact space can be defined as the smallest integer $n$
with the property that any commutative diagram
$$
\xymatrix{ A \ar_{i}[d] \ar^{g}[r]&   K^{{(n)}} \ar_{i}[d]  \\
X \ar^{f}[r]&   K 
}
$$
has a filler $h$
$$
\xymatrix{ A \ar_{i}[d] \ar^{g}[r]&   K^{{(n)}} \ar_{i}[d]  \\
X \ar^{h}[ur]\ar^{f}[r]&   K 
}
$$

Here $A$ is any closed subset of $X$, $K$ is any simplicial complex with the metric topology,
$K^{(n)}$ is the $n$-skeleton of $K$, and $i\colon A\to X$, $i\colon K^{(n)}\to K$ are inclusions.
By saying $h$ is a {\bf filler} we mean $h\vert A=g$ and, since we cannot insist on
$i\circ h=f$, we require $h(x)\in\Delta$ whenever $f(x)\in\Delta$ for any simplex $\Delta$
of $K$.

In \cite{CDV1}, a generalization of the above result was announced for the coarse category of metric
spaces. However, Kevin Zhang, a PHD student in Fudan University of China, found a gap in that paper.
Therefore, the goal of the present paper is not only to provide a proof but to generalize the result
even further, namely to the category of arbitrary coarse spaces. This is done by demonstrating existence
of useful partitions of unity for point-finite covers of coarse spaces (see Section \ref{PartitionsOfUnity}).

In our work we will not use the original 
description of the coarse category of J.Roe. Instead, we will rely on the alternative description
provided in \cite{DH} that is more useful in the context of asymptotic dimension.

The first issue is to find the analog of continuous maps $f\colon X\to K$ from $X$ to a simplicial complex $K$.

As seen in \cite{Dyd} the optimal way to define paracompact spaces $X$
is as follows: for each open cover $\mathcal{U}$ of $X$ there is
a simplicial complex $K$ and a continuous map $f\colon X\to K$ such that
the family $\{f^{-1}(st(v))\}_{v\in K^{(0)}}$ refines $\mathcal{U}$.

In \cite{CDV1} the continuous functions $f:X\to K$ were replaced by
$(\lambda,C)$-Lipschitz functions and the analog of paracompact spaces in coarse geometry
was defined as follows:

\begin{Def}\label{ls-paracompactDef}\cite{CDV1}
A metric space $X$ is {\bf large scale paracompact} (ls-paracompact for short)
if for each uniformly bounded cover $\mathcal{U}$ of $X$ and for all $\lambda,C > 0$
there is a $(\lambda,C)$-Lipschitz function $f\colon X\to K$
such that $\mathcal{V}:=\{f^{-1}(st(v))\}_{v\in K^{(0)}}$ is uniformly bounded
and $\mathcal{U}$ refines $\mathcal{V}$.
\end{Def}

To simplify \ref{ls-paracompactDef} the following concept was introduced:

\begin{Def}\cite{CDV1}
Given $\delta > 0$ and a simplicial complex $K$, a function $f\colon X \to K$
is called a {\bf $\delta$-partition of unity} if it is $(\delta,\delta)$-Lipschitz,
$\mathcal{V}:=\{f^{-1}(st(v))\}_{v\in K^{(0)}}$ is uniformly bounded,
and the Lebesgue number of $\mathcal{V}$
is at least $\frac{1}{\delta}$.
\end{Def}

For arbitrary coarse spaces we need different but related concepts.

\begin{Def}
 Given a cover $\UU$ of a set $X$ and given a function $f:X\to M$ from $X$ to a metric space $M$,
 the \textbf{$\UU$-variation} $var_{\UU}(f)$ of $f$ is the supremum of $d(f(x),f(y))$, where $\{x,y\}$ is contained
 a single element of $\UU$.
  \end{Def}

\begin{Def}
 Given a cover $\UU$ of a coarse space $X$, given $\epsilon > 0$, and given a partition of unity $f:X\to K$,
 we say $f$ is a \textbf{$(\UU,\epsilon)$-partition of unity} if the following conditions are satisfied:\\
a. $var_{\UU}(f) <\epsilon$,\\
b. for every $U\in \UU$ there is $v\in K^{(0)}$  such that
 $f_v(y) > 0$ for all $y\in U$. In other words, point-inverses under $f$ of stars of vertices of $K$
 coarsen $\UU$,\\
 c. point-inverses under $f$ of stars of vertices of $K$ form a uniformly bounded cover of $X$.
 \end{Def}

We are grateful to Kevin Zhang for pointing out a gap in the paper \cite{CDV1}.

\section{Coarsening and shrinking of covers}

In this section we construct shrinking of a coarsening of a cover $\UU$ of $X$
that has multiplicity at most that of $\UU$ and is a coarsening of $\UU$. That
allows us to create useful partitions of unity.

\begin{Def}\label{ShrinkingCoarseningDef}
Given a cover $\UU$ of a set $X$, its \textbf{coarsening} is any cover $\VV$
such that $\UU$ is a refinement of $\VV$.

A \textbf{shrinking} of a cover $\VV=\{V_s\}_{s\in S}$ of $X$ 
is a cover $\WW=\{W_s\}_{s\in S}$ of $X$ such that 
$W_s\subset V_s$ for each $s\in S$.
\end{Def}

\begin{Def}\label{MultiplicityDef}
Given a cover $\UU$ of a set $X$, the \textbf{multiplicity} $m_{\UU}(x)$
of $\UU$ at $x$ is the number of elements of $\UU$ containing $x$.
In case of infinitely many elements of $\UU$ containing $x$ we simply denote it by
$m_{\UU}(x)=\infty$.

A cover $\UU$ is \textbf{point-finite} if $m_{\UU}(x) < \infty$ for all $x\in X$.
\end{Def}

\begin{Lem}\label{ShrinkingWithMultiplicity}
 If $\UU$ is a cover of a set $X$,
 then for every coarsening $\VV=\{V_s\}_{s\in S}$ of $\UU$ there is a shrinking
 $\WW=\{W_s\}_{s\in S}$ of $\VV$ that is a coarsening of $\UU$ such that 
 $m_{\WW}(x)\leq m_{\UU}(x)$ for all $x\in X$.
Moreover, if $x\in V_s\in \VV$ and $m_{\VV}(x)\leq m_{\UU}(x)$, then $x\in W_s$.
\end{Lem}
\dokaz Well-order $S$ and define $W'_s$ as the union of all $U\in \UU$ such that
$s$ is the smallest element of $S$ for which $U\subset V_s$.
$W_s$ is the union of $W'_s$ and of all $x\in V_s$ satisfying $m_{\VV}(x)\leq m_{\UU}(x)$.

Clearly, $\{W_s\}_{s\in S}$ is a shrinking of $\VV$ and a coarsening of $\UU$.
Also, any $x\in X$ that belongs to at least $m_{\UU}(x)+1$ elements of $\{W_s\}_{s\in S}$ must belong
to at least $m_{\UU}(x)+1$ elements of $\UU$, hence $m_{\WW}(x)\leq m_{\UU}(x)$.
\edokaz

\section{Partitions of unity}\label{PartitionsOfUnity}

In this section we construct a partition of unity for every refinement
of a cover of finite multiplicity.

Given a set $S$ of vertices by $\Delta(S)$ we mean the \textbf{full complex}
over $S$: the set of functions $f\colon S\to [0,1]$ with finite support
such that $\sum\limits_{s\in S}f(s)=1$. $\Delta(S)$ is a subset of $l^{1}(S)$,
the space of all functions $f\colon S\to \mathbb{R}$ such that the $l^{1}$-norm $
\| f\|_1=\sum\limits_{s\in S}|f(s)|$ 
of $f$ is finite. $\Delta(S)$ inherits the resulting metric from $l_{1}(S)$.

By a \textbf{simplicial complex} $K$ we mean a subcomplex of $\Delta(S)$ for some set
$S$ ($S$ could be larger than the set of vertices $K^{(0)}$ of $K$).

Any function $g\colon X\to K$ from a space $X$ to a simplicial complex
$K$ can be viewed as a point-finite \textbf{partition of unity} $\{g_v\}_{v\in K^{(0)}}$,
where $g_v(x):=f(x)(v)$.

Given a vertex $v\in K^{(0)}$ by the {\bf star} $st(v)$ of $v$ in $K$
we mean all $f\in K$ such that $f(v) > 0$. Geometrically, it is the union of interiors
of all simplices of $K$ containing $v$.

\begin{Def}
 Given a cover $\UU$ of a set $X$, given $x\in X$, and given $V\subset X$
 we define the \textbf{index $i_{\UU}(x,V)$ of $x$ in $V$ with respect to $\UU$}
 as the smallest integer $k\ge 0$ such that there is a chain of points $x_0=x, x_1,\ldots, x_k$
 with $x_k\notin V$ and $\{x_i,x_{i+1}\}$ belonging to an element of $\UU$ for all $i < k$.
 If such a chain does not exist, we put $i_{\UU}(x,V)=\infty$.
\end{Def}

If the multiplicity function $m_{\VV}$ of a cover $\VV=\{V_s\}_{s\in S}$ is finite at each point, then $\VV$ has a natural partition of unity $\phi_{\UU}^{\VV}$ associated to it via $\UU$:
$$(\phi_{\UU}^{\VV})_s(x)=\frac{i_{\UU}(x,V_s)}{\sum\limits_{t\in S}i_{\UU}(x,V(t))}.$$
In case there are indices $t\in S$ such that $i_{\UU}(x,V(t))=\infty$, we count the number of such indices, say
there is $k$ of them, and we put
$(\phi_{\UU}^{\VV})_s(x)=1/k$ if $i_{\UU}(x,V_s)=\infty$ and $(\phi_{\UU}^{\VV})_s(x)=0$ if $i_{\UU}(x,V_s)< \infty$.

That partition of unity can be considered as a {\bf barycentric map}
$\phi_{\UU}^{\VV}:X\to\NN(\VV)$ from $X$ to the {\bf nerve} of $\VV$. Recall $\NN(\VV)$
is a simplicial complex with vertices belonging to $\VV$
and $\{V_1,\ldots,V_k\}$ is a simplex in $\NN(\VV)$
if and only if $\bigcap\limits_{i=1}^kV_i\ne\emptyset$.

We will be mostly interested in the situation where $\UU$ is a refinement of $\VV$.

\begin{Lem}\label{UVPartitions}
 If $\UU$ is a refinement of $\VV$ and the multiplicity function $m_{\VV}$ of $\VV$ is finite at each point, then
 point-inverses under $\phi_{\UU}^{\VV}$ of stars of vertices of $\NN(\VV)$ coarsen $\UU$.
\end{Lem}
\dokaz Suppose $x\in U\in \UU$. If there is $V_s\in \VV$ such that
$i_{\UU}(x,V_s)=\infty$, then $U\subset V_s$ as otherwise there is $x_1\in U\setminus V_s$
and $i_{\UU}(x,V_s)=1$, a contradiction. For the same reason $i_{\UU}(y,V_s)=\infty$
for all $y\in U$ resulting in $(\phi_{\UU}^{\VV})_s(y) > 0$ and
$y\in (\phi_{\UU}^{\VV})^{-1}(st(V_s))$.

If $i_{\UU}(x,V_s) < \infty$ for all $V_s\in \VV$, we pick $V_s$ such that
$U\subset V_s$. Now $i_{\UU}(x,V_s)\ne 0$ and $(\phi_{\UU}^{\VV})_s(x) > 0$.
Therefore, $x\in (\phi_{\UU}^{\VV})^{-1}(st(V_s))$.
\edokaz

\begin{Lem}\label{MainVarLemma}
 Suppose $\UU$ is a cover of a set $X$, $p:X\to [0,\infty)$, and $q:X\to [m,\infty)$ for some $m > 0$.
 If $p\leq q$, $var_{\UU}(p)\leq 1$, and $var_{\UU}(q)\leq n$, then $var_{\UU}(p/q)\leq \frac{n+1}{m}$.
\end{Lem}
\dokaz
Suppose $x,y\in U\in \UU$. Then $|p(x)-p(y)|\leq 1$, $|q(x)-q(y)|\leq n$, and
$$\left|\frac{p(x)}{q(x)}-\frac{p(y)}{q(y)}\right|=\left|\frac{p(x)\cdot q(y)-p(y)\cdot q(x)}{q(x)\cdot q(y)}\right|=$$
$$\left|\frac{p(x)\cdot (q(y)-q(x))+(p(x)-p(y))\cdot q(x)}{q(x)\cdot q(y)}\right|
\leq \frac{p(x)\cdot n+ q(x)}{q(x)\cdot q(y)}\leq \frac{n}{m}+\frac{1}{m}=\frac{n+1}{m}.$$
\edokaz 

\begin{Def}
 Given a cover $\UU$ of a set $X$ and given $A\subset X$, by $st(A,\UU)$
 we mean the union of all $U\in \UU$ intersecting $A$.
 
 Given two covers $\VV$ and $\UU$ of a set $X$, by $st(\VV,\UU)$ we mean the cover
 $st(V,\UU)$, $V\in \VV$.
 
 $st^k(\UU)$ is defined inductively as follows:\\
 1. $st^0(\UU)=\UU$,\\
 2. $st^{k+1}(\UU)=st(\UU,st^k(\UU))$.
\end{Def}
 
\begin{Cor}\label{MainCorPUs}
Suppose $\VV$ is a uniformly bounded cover of a coarse space $X$.
 If $st^k(\UU)$ refines $\VV$ for some $k\ge 1$ and the multiplicity of $\VV$ is at most $n+1$,
 then the partition $\phi_{\UU}^{\VV}:X\to \NN(\VV)$ is
 a $(\UU,\frac{(2n+2)^2}{k})$-partition of unity.
\end{Cor}
\dokaz Use \ref{UVPartitions} to see that point-inverses of stars of vertices of $\NN(\VV)$
coarsen $\UU$ and refine $\VV$.
Assume $\VV=\{V_s\}_{s\in S}$. If $x$ is a point such that
$\sum\limits_{t\in S}i_{\UU}(x,V(t))=\infty$, then for any $y$ belonging to the same
element $U$ of $\UU$ either both $i_{\UU}(x,V(t))$, $i_{\UU}(y,V(t))$ are finite
or both are infinite. Therefore $\phi_{\UU}^{\VV}(x)=\phi_{\UU}^{\VV}(y)$ and $|\phi_{\UU}^{\VV}(x)-\phi_{\UU}^{\VV}(y)|=0$.

Assume $\sum\limits_{t\in S}i_{\UU}(x,V(t))< \infty$ and $x,y\in U\in \UU$.
Notice there is a subset $T$ of $S$ containing at most $2n+1$ elements such that for
$t\in S\setminus T$ the indices $i_{\UU}(x,V(t))$, $i_{\UU}(y,V(t))$ are $0$.
Therefore $var_{\UU}(\sum\limits_{t\in S}i_{\UU}(x,V(t)))\leq 2n+1$.
Since one of the indices is at least $k$, $\sum\limits_{t\in S}i_{\UU}(x,V(t))\ge k$.
Applying \ref{MainVarLemma} one gets that the $\UU$-variation
of each $(\phi_{\UU}^{\VV})_s(x)=\frac{i_{\UU}(x,V_s)}{\sum\limits_{t\in S}i_{\UU}(x,V(t))}$ is at most
$\frac{2n+2}{k}$. As there are at most $2n+1$ functions that are relevant
for $\phi_{\UU}^{\VV}(x)$ and $\phi_{\UU}^{\VV}(y)$, $|\phi_{\UU}^{\VV}(x)-\phi_{\UU}^{\VV}(y)|\leq \frac{(2n+2)^2}{k}$.
\edokaz

\begin{Lem}\label{ComparisonOfPUs}
 Suppose $X$ is a metric space, $2 > \delta > 0$, and $\UU$ is the cover of $X$ 
 by balls of radius $\frac{1}{\delta}$. Every $\delta^2/4$-partition of unity on $X$
 is a $(\UU,\delta)$-partition of unity on $X$ and every $(\UU,\delta)$-partition of unity on $X$
 is a $2\delta$-partition of unity on $X$.
\end{Lem}
\dokaz Suppose $f:X\to K$ is a $\delta^2/4$-partition of unity. Since $\delta^2/4 < \delta$,
point-inverses of stars of vertices under $f$ refine $\UU$.
If $x,y\in U\in \UU$, then $d(x,y) < 2/\delta$ and
$|f(x)-f(y)|\leq \delta^2\cdot d(x,y)/4+\delta^2/4 < \delta/2+\delta^2/4 < \delta$.

Suppose $\UU$ consists of sets of diameter at most $\frac{2}{\delta}$, is of Lebesgue number
at least $\frac{1}{\delta}$, and $f:X\to K$ is a $(\UU,\delta)$-partition of unity. 
If $d(x,y)\ge \frac{1}{\delta}$, then $|f(x)-f(y)|\leq 2 < 2\delta\cdot \frac{1}{\delta}+2\delta \leq 2\delta\cdot d(x,y)+2\delta $.
If $d(x,y)< \frac{1}{\delta}$, then there is $U\in \UU$ containing both $x$ and $y$
and $|f(x)-f(y)| < \delta \leq 2\delta\cdot d(x,y)+2\delta$.
Thus $f$ is a $2\delta$-partition of unity.
\edokaz

\section{Asymptotic dimension}

\begin{Def}
 A coarse space $X$ has \textbf{asymptotic dimension} $\asdim(X)$ at most $n$ if for every uniformly bounded
 cover $\UU$ of $X$ there is a uniformly bounded cover $\VV$ of $X$ such that
 every element $U$ of $\UU$ intersects at most $n+1$ elements of $\VV$. 
\end{Def}

\begin{Thm}\label{CharOfAsdimViaMapsToSimplicialComplexes}
 The following conditions are equivalent for any coarse space $X$ and any integer $n\ge 0$:\\
 a. $\asdim(X)\leq n$.\\
 b. for every $\epsilon > 0$ and every uniformly bounded cover $\UU$ of $X$
 there is a $(\UU,\epsilon)$-partition of unity $f:X\to K^{(n)}$.\\
 c. for every uniformly bounded cover $\UU$ of $X$
 there is a $(\UU,\infty)$-partition of unity $f:X\to K^{(n)}$.
\end{Thm}
\dokaz
a)$\implies$b) follows from \ref{MainCorPUs}. Indeed, we choose
$k > 1$ so that $\frac{(2n+2)^2}{k} < \epsilon$, and then we choose a uniformly bounded cover
$\VV$ of multiplicity at most $n+1$ that coarsens $st^k(\UU)$. The partition of unity $\phi_{\UU}^{\VV}$ for that cover is what we need. To obtain $\VV$ we first choose a uniformly bounded cover $\WW$ with the property
that every element of $st^{k+1}(\UU)$ intersects at most $n+1$ elements of $\WW$.
Now, $\VV:=st(\WW,st^k(\UU))$ works: if a point $x\in X$ belongs to
$n+2$ elements of $\VV$, then $st(x,st^k(\UU))$ intersects $n+2$ elements of $\WW$,
a contradiction.
\\
b)$\implies$c) is obvious.\\
c)$\implies$a). Given a uniformly bounded cover $\UU$ of $X$ pick a $(st^2(\UU),\infty)$-partition of unity $f:X\to K^{(n)}$.
Obviously, point-inverses under $f$ of stars of vertices of $K$ form a uniformly bounded cover of $X$
that coarsens $st^2(\UU)$ and is of multiplicity at most $n+1$.
Remove from each $f^{-1}(st(v))$ the union of all $U\in \UU$ intersecting the complement of $f^{-1}(st(v))$.
The new family $\VV$ is a coarsening of $\UU$.
Suppose there is $U\in \UU$ intersecting at least $n+2$ elements of $\VV$.
That implies $U$ is contained in $n+2$ different point-inverses $f^{-1}(st(v))$, a contradiction.
\edokaz

\begin{Thm}\label{TheoremANew}
Suppose $X$ is a coarse space.
 If $X$ is of asymptotic dimension at most $n\ge 0$, then for every $\epsilon > 0$
 and every uniformly bounded cover $\UU$ of $X$
 there is $\delta > 0$ and a uniformly bounded cover $\VV$ of $X$ such that any commutative diagram
 $$
\xymatrix{ A \ar_{i}[d] \ar^{g}[r]&   K^{{(n)}} \ar_{i}[d]  \\
X \ar^{f}[r]&   K 
}
$$
where $f$ is a $(\VV,\delta)$-partition of unity
has a filler $h$
$$
\xymatrix{ A \ar_{i}[d] \ar^{g}[r]&   K^{{(n)}} \ar_{i}[d]  \\
X \ar^{h}[ur]\ar^{f}[r]&   K 
}
$$
that is a $(\UU,\epsilon)$-partition of unity.
\end{Thm}
\dokaz 
Choose $k > m$ and $\delta$ such that each term of the sum
$$(2m+1)\delta+3/m+2(2n+2)^2/k+3/m$$
 is smaller than $\min(\epsilon,\frac{1}{n+1})/4$.
$\VV$ is chosen so that it coarsens $st^k(\UU)$ and is of multiplicity
at most $n+1$.

Define $\alpha(x)$ as $\frac{i_{\UU}(x,st^m(A,\UU))}{i_{\UU}(x,st^m(A,\UU))+i_{\UU}(x,X\setminus A)}$
if both indices are finite. Define $\alpha(x)=1/2$ if both indices are infinite.
If $i_{\UU}(x,st^m(A,\UU))=\infty$ and $i_{\UU}(x,X\setminus A) < \infty$, we put $\alpha(x)=1$.
If $i_{\UU}(x,st^m(A,\UU)) < \infty$ and $i_{\UU}(x,X\setminus A) = \infty$, we put $\alpha(x)=0$.

 Notice $var_{U}(\alpha)\leq 3/m$ for points where both indices are finite
by \ref{MainVarLemma}. Indeed, $i_{\UU}(x,st^m(A,\UU))+i_{\UU}(x,X\setminus A)\ge m$ for each $x\in st^m(A,\UU)$. 
Also, $i_{\UU}(x,st^m(A,\UU)) > m$ for each $x\notin st^m(A,\UU)$. 

In case of points $x,y\in U\in \UU$ such that at least one index is infinite, we observe that
the corresponding indices are either both finite or both infinite. Therefore $\alpha(x)=\alpha(y)$
and $var_{\UU}(\alpha)\leq 3/m$ holds in full generality.

Given $f:X\to K$ such that $f(A)\subset K^{(n)}$ we find for each $x\in st^m(A,\UU)$ a point
$c(x)\in A$ that can be reached by the shortest $\UU$-chain. For $x\in A$ we put $c(x)=x$.
Since that chain has at most $m$ links,
$|f(x)-f(c(x))|\leq m\cdot \delta < (\frac{1}{n+1})/8$. In particular, there are common vertices $v$
in the carriers of $f(x)$ and $f(c(x))$ (i.e. $f_v(x), f_v(c(x)) > 0$).

We define $g:st^m(A,\UU)\to K^{(n)}$ by adding all
coefficients of $f(x)$ that are not in the carrier of $f(c(x))$ to one vertex of the intersection of the two carriers.
Thus $|g(x)-f(x)|\leq m\cdot \delta$ and $var_{\VV}(g)\leq (2m+1)\cdot \delta$.

We shrink $f^{-1}(st(v))$, $v\in K^{(0)}$, as in \ref{ShrinkingWithMultiplicity}  to a cover $\WW=\{W_v\}$ coarsening $\VV$
of multiplicity at most $n+1$ (see \ref{ShrinkingWithMultiplicity}). The partition of unity $\phi:=\phi_{\UU}^{\WW}:X\to K$ for that cover
has $\UU$-variation at most $(2n+2)^2/k$ by \ref{MainCorPUs}.

Finally, we create $h(x)=\alpha(x)\cdot g(x)+(1-\alpha(x))\cdot \phi(x)$. To estimate its $\UU$-variation
notice, given $x,y$ belonging to $U\in \UU$, that
$$h(x)-h(y)=\alpha(x)(g(x)-g(y))+(\alpha(x)-\alpha(y))\cdot g(y)+(\phi(y)-\phi(x))+$$
$$\alpha(x)(\phi(y)-\phi(x))+(\alpha(y)-\alpha(x))\cdot \phi(y).$$
Consequently,
$$|h(x)-h(y)|\leq (2m+1)\delta+3/m+2(2n+2)^2/k+3/m < \min(\epsilon,\frac{1}{n+1}).$$
In particular, if $x\in U\in \VV$ and $h_v(x)\ge 1/(n+1)$ (such a vertex $v$ of $K$ exists),
then for each $y\in U$ we have $h_v(y) > 0$. That shows point inverses of stars of vertices of $K$
under $h$ coarsen $\UU$ and $h$ is a $(\UU,\epsilon)$-partition of unity.
\edokaz

We are now able to provide a proof of the result announced in \cite{CDV1}.

\begin{Cor}\label{TheoremA}\cite{CDV1}
 If $X$ is of asymptotic dimension at most $n\ge 0$, then for any $\epsilon > 0$
 there is $\delta > 0$ such that any commutative diagram
 $$
\xymatrix{ A \ar_{i}[d] \ar^{g}[r]&   K^{{(n)}} \ar_{i}[d]  \\
X \ar^{f}[r]&   K 
}
$$
where $f$ is a $\delta$-partition of unity
has a filler $h$
$$
\xymatrix{ A \ar_{i}[d] \ar^{g}[r]&   K^{{(n)}} \ar_{i}[d]  \\
X \ar^{h}[ur]\ar^{f}[r]&   K 
}
$$
that is an $\epsilon$-partition of unity.
\end{Cor}
\dokaz
Apply \ref{ComparisonOfPUs} and \ref{TheoremANew}.
\edokaz

\section{Large scale paracompactness}
 
 Large scale paracompact spaces were introduced in \cite{CDV2}
 as a generalization of metric spaces with Property A (see \cite{NY} and 
\cite{Willett}). In this section we show that the class of large scale paracompact spaces
plays the same role in the coarse category (when it comes to asymptotic dimension) as the class of paracompact spaces.

\begin{Def}
 A coarse space $X$ is \textbf{large scale paracompact} if for every uniformly bounded cover $\UU$ of $X$
 and every $\epsilon > 0$ there is a $(\UU,\epsilon)$-partition of unity $f:X\to K$.
\end{Def}

\begin{Cor}\label{TheoremBNew}
Suppose $X$ is a large scale paracompact space and $n\ge 0$. If for every $\epsilon > 0$
 and every uniformly bounded cover $\UU$ of $X$
 there is $\delta > 0$ and a uniformly bounded cover $\VV$ of $X$ such that any commutative diagram
 $$
\xymatrix{ A \ar_{i}[d] \ar^{g}[r]&   K^{{(n)}} \ar_{i}[d]  \\
X \ar^{f}[r]&   K 
}
$$
where $f$ is a $(\VV,\delta)$-partition of unity
has a filler $h$
$$
\xymatrix{ A \ar_{i}[d] \ar^{g}[r]&   K^{{(n)}} \ar_{i}[d]  \\
X \ar^{h}[ur]\ar^{f}[r]&   K 
}
$$
that is a $(\UU,\epsilon)$-partition of unity, then $\asdim(X)\leq n$.
\end{Cor}
\dokaz
Suppose $\UU$ is a uniformly bounded cover of $X$.
Choose a uniformly bounded cover $\VV$ and $\delta > 0$ for $\epsilon=1$.
As $X$ is large scale paracompact, there is a $(\VV,\delta)$-partition of unity
$f\colon X\to K$.
Let $h\colon X\to K^{(n)}$ be a $(\UU,\epsilon)$-partition of unity and a filler
of $f$. Apply \ref{CharOfAsdimViaMapsToSimplicialComplexes}.
\edokaz

Here is another result announced in \cite{CDV1}.
\begin{Cor}\label{TheoremB}\cite{CDV1}
 Suppose $n\ge 0$ and for any $\epsilon > 0$
 there is $\delta > 0$ such that any commutative diagram
 $$
\xymatrix{ A \ar_{i}[d] \ar^{g}[r]&   K^{{(n)}} \ar_{i}[d]  \\
X \ar^{f}[r]&   K 
}
$$
where $f$ is a $\delta$-partition of unity
has a filler $h$
$$
\xymatrix{ A \ar_{i}[d] \ar^{g}[r]&   K^{{(n)}} \ar_{i}[d]  \\
X \ar^{h}[ur]\ar^{f}[r]&   K 
}
$$
that is an $\epsilon$-partition of unity.
If $X$ is large scale paracompact, then its asymptotic dimension is at most $n$.
\end{Cor}
\dokaz
Apply \ref{ComparisonOfPUs} and \ref{TheoremBNew}.
\edokaz

\end{document}